\documentclass [twoside,10pt]{article}
\usepackage{graphicx}
\usepackage{amssymb}
\usepackage{amsthm}

\newtheorem{pro}{Proposition}

\begin{document}

\baselineskip 18.5pt

\title{\vspace{-1cm}{\bf A note on some generalizations of Monge's theorem}}

\pagestyle{myheadings} \markboth{A note on generalizations Monge's theorem}{Marek Lassak}

\author{Marek Lassak}

\date{} 

\maketitle

\vskip0.2cm

\noindent
{\bf Abstract.}
We generalize Monge's theorem for $n+1$ pairwise homothetic sets (in particular convex bodies) in $E^n$ in place of 
three disks in $E^2$. 
We also present a version for $n+1$ independent points of $E^n$.
It also includes the reverse statement.
Moreover, we give an analogon of Monge's theorem for the $n$-dimensional sphere and hyperboloid model of the hyperbolic space.

\vskip0.2cm
\noindent
{\bf Mathematical Subject Classification (2010).} 52A20, 52A21, 52A55. 

\vskip0.2cm
\noindent
{\bf Keywords.} Monge's theorem, Menelaus' theorem, Euclidean space, sphere, hyperbolic space

\section{Introduction}

For any two disjoint disks in a plane, an external tangent is a line that is tangent to both disks but does not pass between them. There are two such external tangent lines for any two circles. 
Each such pair of external tangents for disks of different size has a unique intersection point. 
The classic Monge's theorem states that for three such pairwise disjoint circles of different size the three intersection points of the external tangent lines given by the three pairs of circles always lie in a straight line.
For instance, see the book \cite{G} by Gardner.
Walker \cite{W} generalized Monge's Theorem for $n+1$ balls in the Euclidean $n$-space $E^n$ and a hyperplane in place of the above straight line.

First we present a version of Monge's theorem for $n+1$ linearly independent points in place of the balls.
This version includes also the reverse claim.
Next we give a generalization of Monge's theorem for $n+1$ pairwise homothetic bounded sets (not obligatory disjoint) with homothety ratios over $1$ in $E^n$ in place of the $n+1$ balls in $E^n$.
A good visualization is obtained by taking convex bodies in place of our bounded sets.

Moreover, we give spherical and hyperbolic $n$-dimensional analogons of Monge's theorem.
Our proof applies the $n$-dimensional Menelaus' theorem for the $n$-dimensional sphere $S^n$ and hyperbolic space $H^n$ by Ushijma \cite{U} .
A question remains about possible analogs of Monge's theorem in axiomatic geometry, where Guggenheimer \cite{G} considers the Menelau's theorem.

\section
{Monge's theorem for a wide class of sets in $E^n$}

Let us start with a proposition presenting a version of Monge's theorem for $n+1$ points in $E^n$ instead of balls.

Recall that a set $A$ of points in $E^n$ is said to be {\it independent} if the affine span of any proper subset of $A$ is a proper subset of the span of $A$.

\vskip0.15cm
\begin{pro} 
{\it Let the set of points $a_1, a_2, \dots , a_{n + 1} \in E^n$ be independent.
Consider the straight line $L_{ij}$ containing $a_ia_j$ and a point $b_{ij} \in L_{ij}$ different from $a_i$ and $a_j$ for $i, j \in \{1, \dots , n + 1 \}$.
For every $i<j$ denote by $\lambda_{ij}$ the ratio of homothety with center $b_{ij}$ which transforms $a_j$ into $a_i$.
We claim that the $n(n + 1) / 2$ points $b_{ij}$ belong to one hyperplane if and only if
$\lambda_{ij}^{-1}\lambda_{ik}\lambda_{jk}^{-1} = 1$ for every $i, j \in \{1, \dots , n+1\}$.}
\end{pro} 

\begin{proof}
Let us apply the variant of Theorem 2 of the paper \cite{B} by Buba-Brzozowa in which we take into account just lengths instead of the oriented lengths. 
Her $n$-dimensional generalization of the classic  Menelaus' theorem says that points $b_{ij}$, where $i, j = 1, \dots , n + 1$ and $i<j$, belong to one hyperplane of $E^n$ if and only if 

\vskip-0.4cm

$$\frac{|a_ib_{ij}|}{|b_{ij}a_j|}
\cdot
\frac{|a_jb_{jk}|}{|b_{jk}a_k|}
\cdot
\frac{|a_kb_{ik}|}{|{b_{ik}a_i}|}
= 1.$$

Since $\frac{|a_ib_{ij}|}{|b_{ij}a_j|} = \lambda_{ij}$ for $i, j \in \{1, \dots , n + 1\}$ and $i<j$, we obtain that all our points $b_{ij}$ are in one hyperplane if and only if $\lambda_{ij}^{-1}\lambda_{ik}\lambda_{jk}^{-1} = 1$ for every $i, j \in \{ 1, \dots , n + 1 \}$, which is our thesis.  
\end{proof}

\vskip0.15cm
Clearly, if we agree that the different points $a_1, \dots, a_{n + 1}$ in Proposition 1 are dependent, then the ``if" part trivially holds true.

\vskip0.45cm
\noindent
{\bf Theorem.}
{\it Assume that for sets $C_1, \dots , C_{n + 1} \subset E^n$ and for every  $i, j \in \{1, \dots , n + 1\}$ with $i<j$ there are unique homotheties $h_{ij}$ of ratios over $1$ such that $h_{ij}(C_j) = C_i$.
Then the $n(n + 1) / 2$ centers of these homotheties are in one hyperplane.}

\begin{proof}
Consider the $n(n + 1) / 2$ homotheties $h_{ij}$ such that $h_{ij}(C_j) = C_i$, where $i, j \in \{1, \dots , n + 1 \}$ and $i<j$ (see Figure for the special case of three convex bodies in $E^2$). 
Of course, for every three homotheties $h_{pq}, h_{pr}, h_{qr}$, where $p<q<r$, we have $h_{pr}(h_{qr}^{-1}(h_{pq}^{-1}(C_p))) = C_p$.
By the assumed uniqueness of the homotheties $h_{ij}$ there are points $o_1, \dots , o_{n + 1}$ such that $h_{ij} (o_j) = o_i$ for every $i, j \in \{ 1, \dots , n + 1 \}$ with $i < j$. 
From the ``if" part of Proposition 1 we get our thesis.
\end{proof}

{\ }

\begin{center}

{\ }

\vskip-0.25cm 

\includegraphics[width=4.3in]{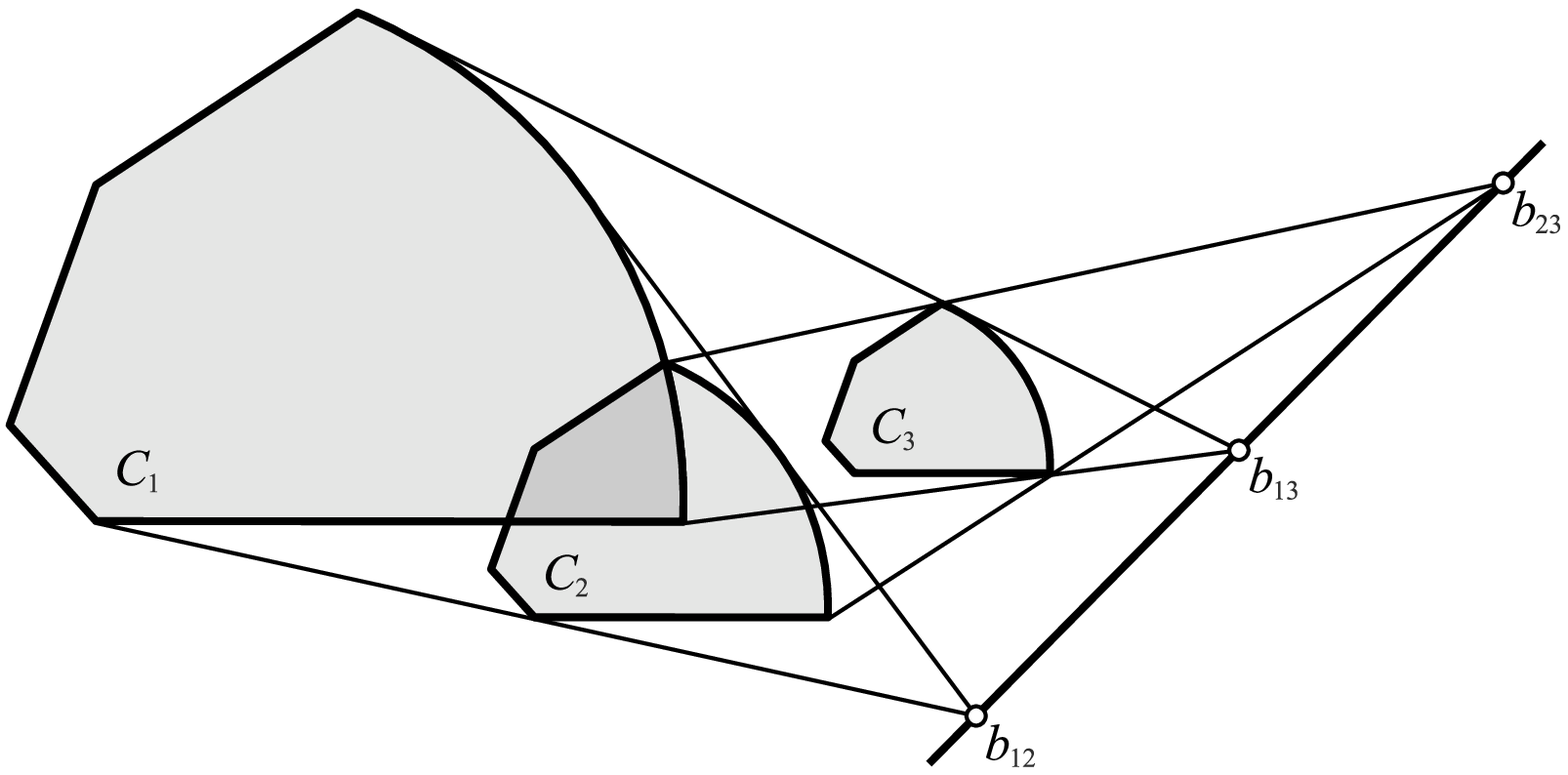} \\ 

\vskip 0.05cm
{Figure. Illustration of Theorem for the case of three convex bodies in $E^2$}  

\end{center}

\vskip0.1cm
If $C_1, \dots , C_{n + 1}$ from Theorem are centrally-symmetric convex bodies, we may rephrase this as follows: {\it Monge's theorem holds for $n + 1$ balls of different sizes of a normed $n$-dimensional space} in place of $n + 1$ balls of Euclidean space.
A particular case is for the two-dimensional $L_p$ spaces considered in the preprint \cite{EG} by Ermi\c s and Geli\c sgen.
The author thanks them for sharing their preprint \cite{EG}, which was mobilizing to show the above theorem.

The assumption that homotheties $h_{ij}$ are unique holds always if the sets are bounded and non-empty. 
It also holds for some unbounded sets.
For instance when in $E^2$ we take $C_i$ as the intersection of half-planes $x\geq 0$, $\ y \geq 1/i$ and $x + y \geq 4 - i$, where $i = 1, 2, 3$.
Then all $b_{ij}$ are different and lie on the line $x = 0$.
If we exchange $y \geq 1/i$ into $ y \geq 0$ here, all $b_{ij}$ coincide and still are on $x = 0$.
For example the assumption does not hold for any family of $n + 1$ translates of a half-space in $E^n$;
the thesis may be not true for some homotheties with ratios over 1 between them.
We let the reader to show that if the latter assumption does not hold, then there are some positive homotheties for which the thesis of Theorem is still true.

\section {Analogons of Monge's theorem in spherical and hyperbolic spaces}

Below by $X^n$ we denote both the $n$-dimensional sphere $S^n$ and the hyperboloid model $H^n$ of the hyperbolic $n$-dimensional space.
By a {\it hyperplane} and a {\it line} of $X^n$ we mean a subset of $X^n$ isometric to $X^{n-1}$ and $X^1$, respectively.
Usually, for $S^n$ they are called a {\it $(n-1)$-dimensional subsphere} and a {\it great circle}, respectively.
By the {\it distance} $|xy|$ of points $x,y \in X^n$ (which are not opposite for $X^n = S^n$) we mean the length of the geodesic joining $x$ and $y$.
By the arc $xz$ between $x$ and $z$ we mean the set of points $y$ such that $|xy| + |yz| = |xz|$. 
Let $c, p, r \in X^n$ be points such that $p \in cr$ and $c \not = p$,  or $r \in cp$ and $c \not = r$.
If $|cr| = \lambda \cdot |cp|$, then we say that $r$ is the {\it image} of $p$ under the {\it $X^n$-homothety} with center $c$ and ratio $\lambda$.
Clearly $\lambda > 0$.

We call a set $A \subset X^n$  
(embedded in $E^{n + 1}$) to be {\it independent} if the set $A \cup \{o\}$, where $o$ is the origin of $E^{n + 1}$, is independent in  $E^{n + 1}$. 

The Menelau's theorem on $S^2$ is recalled in Proposition 66 of the book \cite{RP} by Rashed and Papadopoulos and its variant for $H^2$ is presented by Smarandache and Barbu \cite{SB}. 
Recently their generalizations for $S^n$ and $H^n$ are giver in Theorem 4 of Ushijima \cite{U}.
From this result, similarly to the proof of our Proposition 1 for $E^n$, we get the following Proposition 2 on the analogous variant of Monge's theorem for $n + 1$ points on $S^n$ and $H^n$.
Here by $\lambda_{ij}$ we mean $\frac{\sin |a_ib_{ij}|}{\sin |b_{ij}a_j|}$ for $S^n$ and $\frac{\sinh |a_ib_{ij}|}{\sinh |b_{ij}a_j|}$ for $H^n$.

\begin{pro} 
{\it Let $a_1, a_2 \dots , a_{n + 1}$ be a set of independent points of $X^n$. 
Consider the line $L_{ij}$ containing $a_ia_j$. 
Denote by $b_{ij}$ a point different from $a_i$ and $a_j$ in $L_{ij}$ such that $a_j \in a_ib_{ij}$ for $i, j \in \{ 1, \dots n + 1 \}$.
For every $i<j$ denote by $\lambda_{ij}$ the ratio of the $X^n$-homothety with center $b_{ij}$ such that $a_i$ is the image of $a_j$.
Then the $n(n +1 )/2$ points $b_{ij}$ are in one hyperplane of $H^n$ if and only if
$\lambda_{ij}^{-1}\lambda_{ik}\lambda_{jk}^{-1} = 1$ for every $i, j \in \{1, \dots , n + 1\}$.} 
\end{pro}

\baselineskip 12pt

\vskip0.1cm
\noindent
Marek Lassak

\noindent
University of Science and Technology

\noindent
Bydgoszcz 85-798, al. Kaliskiego 7, Poland

\noindent
e-mail: lassak@utp.edu.pl

\end{document}